# Two-Person Additively-Separable Sum Games

By

Somdeb Lahiri

(somdeb.lahiri@gmail.com)

ORCID: https://orcid.org/0000-0002-5247-3497

**Abstract**

We consider a sub-class of bi-matrix games which we refer to as two-person (hereafter referred to as two-player) additively-separable sum (TPASS) games, where the sum of the pay-offs of the two players is additively separable. The row player's pay-off at each pair of pure strategies, is the sum of two numbers, the first of which may be dependent on the pure strategy chosen by the column player and the second being independent of the pure strategy chosen by the column player. The column player's pay-off at each pair of pure strategies, is also the sum of two numbers, the first of which may be dependent on the pure strategy chosen by the row player and the second being independent of the pure strategy chosen by the row player. The sum of the inter-dependent components of the pay-offs of the two players is assumed to be zero. We prove the existence of equilibrium for such games and show that the set of equilibria for such games is the projection on the set of strategy pairs of the solutions of a pair of linear programming problem that are dual to each other. This result is a generalization of the corresponding and well-known result for two-person zero-sum games. We also show that a (randomized or mixed) strategy pair is an equilibrium of the game if and only if there exist two other real numbers such that the three together solve a certain linear programming problem. In order to prove this result, we need to appeal to the existence of an equilibrium for the TPASS game. The technology we use to prove our results, consists of the duality theorems and the complementary slackness theorem of linear programming.

**1.** Consider an ordered triplet $(A, \pi, \rho)$ where for some positive integers m and n, A is an $m \times n$ real valued matrix, $\pi$ is an m-dimensional real-valued column vector and $\rho$ is an n-dimensional real-valued column vector. For $(i, j) \in \{1, \ldots, m\} \times \{1, \ldots, n\}$, let $\pi_i$ denote the $i^{th}$ coordinate of $\pi$, $\rho_j$ denote the $j^{th}$ coordinate of $\rho$, and let $a_{ij}$ denote the entry at the intersection of the $i^{th}$ row and $j^{th}$ column of A.

There are two players in this game (i.e., interactive decision-making problem)- the row player and the column player.

For $(i, j) \in \{1, \ldots, m\} \times \{1, \ldots, n\}$, if the row player chooses the $i^{th}$ row and the column player chooses column j, then the payoff to the row player is $a_{ij} + \pi_i$ and the payoff to the column player is $- a_{ij} + \rho_j$.

Note that the sum of the pay-offs to the row player and the column player if the row player chooses the $i^{th}$ row and the column player chooses $j^{th}$ column is $\pi_i + \rho_j$.

We refer to the triplet $(A, \pi, \rho)$ as a **two-person additively-separable sum** (TPASS) game.

In case all coordinates of $\pi$ are identical and all coordinates of $\rho$ are identical, a TPASS game reduces to a two-person constant sum (TPCS) game (see chapter 20 of Mote and Madhavan (2016)).

We allow for randomized (mixed) strategies for the row and column players.

For any non-negative integer $\ell$, let $\Delta^{\ell} = \{x \in \mathbb{R}_+^{\ell+1} | \sum_{k=1}^{\ell+1} x_k = 1\}$. We will interpret points in $\Delta^{\ell}$ as $\ell$-dimensional column vectors.

The (randomized or mixed) strategy set for the row player is $\Delta^{m-1}$ and the (randomized or mixed) strategy set for the column player is $\Delta^{n-1}$.

A pair $(p, q) \in \Delta^{m-1} \times \Delta^{n-1}$ is a (**randomized or mixed**) **strategy pair**.

The pay-off function for the row-player is the function $f^R: \Delta^{m-1} \times \Delta^{n-1} \rightarrow \mathbb{R}$ such that for all $(p, q) \in \Delta^{m-1} \times \Delta^{n-1}$, $f^R(p, q) = p^T A q + p^T \pi$.

The pay-off function for the column-player is the function $f^C: \Delta^{m-1} \times \Delta^{n-1} \rightarrow \mathbb{R}$ such that for all $(p, q) \in \Delta^{m-1} \times \Delta^{n-1}$, $f^C(p, q) = -p^T A q + \rho^T q$.

The following concept is available in Nash (1951).

$(p^*, q^*) \in \Delta^{m-1} \times \Delta^{n-1}$ is said to be **an equilibrium** of the TPASS game $(A, \pi)$ if for all $(p, q) \in \Delta^{m-1} \times \Delta^{n-1}$: $f^R(p^*, q^*) \geq f^R(p, q^*)$ and $f^C(p^*, q^*) \geq f^C(p^*, q)$.

Let $R(\pi)$ be the $m \times n$ real matrix, such that for all $i \in \{1, \ldots, m\}$, every entry in the $i^{th}$ row of $R(\pi)$ is $\pi_i$.

Let $C(\rho)$ be the $m \times n$ real matrix, such that for all $j \in \{1, \ldots, n\}$, every entry in the $j^{th}$ column of $C(\rho)$ is $\rho_j$.

Thus, for all $(p, q) \in \Delta^{m-1} \times \Delta^{n-1}$, $p^T \pi = p^T R(\pi) q$ and $\rho^T q = p^T C(\rho) q$.

For any positive integer $\ell$, let $e^{(\ell)}$ denote the $\ell$-dimensional sum column vector, i.e., the $\ell$-dimensional vector, all coordinates of which are 1.

**2.** Before proving the characterization theorem for the set of all equilibria of a TPASS game, let us prove the existence of equilibrium for TPASS games. Our proof, as in the case of TPCS games, uses no more than the duality and complementary slackness theorems for linear programming problems.

**Proposition 1:** Let $(A, \pi, \rho)$ be a TPASS game. Then, it has an equilibrium.

**Proof:** Consider the linear programming problem LP1:

Maximize $\rho^T q - \alpha$, subject to $Aq - \alpha e^{(m)} \leq - \pi$, $e^{(n)T}q = 1$, $q \in \mathbb{R}^n_+$, $\alpha \in \mathbb{R}$.

The dual of LP1 is the linear programming problem Dual-LP1:

Minimize $-\pi^T p + \beta$, subject to $p^T A + \beta e^{(n)T} \geq \rho^T$, $-p^T e^{(m)} = -1$, $p \in \mathbb{R}^m_+$, $\beta \in \mathbb{R}$.

Clearly both LP1 and Dual-LP1 have feasible solutions. Hence, by the strong duality theorem of linear programming, both have optimal solution. Let $q^*$, $\alpha^*$ solve LP1 and $p^*$, $\beta^*$ solve Dual-LP1.

Then, from feasibility and the well know complementary slackness conditions we know that the following conditions are satisfied:

(i) $Aq^* - \alpha^* e^{(m)} \leq - \pi$, $e^{(n)T}q^* = 1$, $q^* \in \mathbb{R}^n_+$, $\alpha^* \in \mathbb{R}$.

(ii) $p^{*T}A + \beta^* e^{(n)T} \geq \rho^T$, $e^{(m)T}p^* = 1$, $p^* \in \mathbb{R}^m_+$, $\beta \in \mathbb{R}$.

(iii) $p^{*T}Aq^* - \alpha^* = -p^{*T}\pi$, $p^{*T}Aq^* + \beta^* = \rho^T q^*$.

From (iii) we get $\alpha^* = p^{*T}Aq^* + p^{*T}\pi$ and $- p^{*T}Aq^* + \rho^T q^* = \beta^*$.

Thus, these two equations along with (i) and (ii) imply $Aq^* + \pi \leq (p^{*T}Aq^* + p^{*T}\pi)e^{(m)}$, $-p^{*T}A + \rho^T \leq (- p^{*T}Aq^* + \rho^T q^*)e^{(n)T}$, $(p^*, q^*) \in \Delta^{m-1} \times \Delta^{n-1}$.

Thus, $(p^*, q^*) \in \Delta^{m-1} \times \Delta^{n-1}$ and for all $(p, q) \in \Delta^{m-1} \times \Delta^{n-1}$: $f^R(p^*, q^*) \geq f^R(p, q^*)$ and $f^C(p^*, q^*) \geq f^C(p^*, q)$.

Thus, $(p^*, q^*)$ is an equilibrium for the TPASS game $(A, \pi, \rho)$. Q.E.D.

It is easy to see that any equilibrium of the TPASS game $(A, \pi, \rho)$ along with appropriate scalars provides solutions for the linear programming problem in the proof of proposition 1 and its dual. Thus, we have the following result.

**Proposition 2:** $(p^*, q^*)$ is an equilibrium for the TPASS game $(A, \pi, \rho)$ if and only if there exist real numbers $\alpha^*$, $\beta^*$ such that $q^*$, $\alpha^*$ solve [Maximize $\rho^T q - \alpha$, subject to $Aq - \alpha e^{(m)} \leq -\pi$, $e^{(n)T}q = 1$, $q \in \mathbb{R}^n_+$, $\alpha \in \mathbb{R}$] and $p^*$, $\beta^*$ solve its dual [Minimize $-\pi^T p + \beta$, subject to $p^T A + \beta e^{(n)T} \geq \rho^T$, $-p^T e^{(m)} = -1$, $p \in \mathbb{R}^m_+$, $\beta \in \mathbb{R}$]. Further, $p^*$, $q^*$, $\alpha^*$, $\beta^*$ satisfies $\alpha^* = p^{*T}Aq^* + p^{*T}\pi$ and $\beta^* = - p^{*T}Aq^* + \rho^T q^*$.

**Proof:** If q, $\alpha$ satisfy [$Aq - \alpha e^{(m)} \leq - \pi$, $e^{(n)T}q = 1$, $q \in \mathbb{R}^n_+$, $\alpha \in \mathbb{R}$] and p, $\beta$ satisfy [$p^T A + \beta e^{(n)T} \geq \rho^T$, $-p^T e^{(m)} = -1$, $p \in \mathbb{R}^m_+$, $\beta \in \mathbb{R}$] then since the two-linear programming problems in the statement of the proposition are dual to each other, it follows that $\rho^T q - \alpha \leq -\pi^T p + \beta$.

**Part 1:** Let $(p^*, q^*)$ be an equilibrium for the TPASS game $(A, \pi, \rho)$. Thus, $Aq^* + \pi \leq (p^{*T}Aq^* + p^*\pi)e^{(m)}$ and $-p^{*T}A + \rho^T \leq (-p^{*T}Aq^* + \rho^T q^*)e^{(n)T}$.

Choose real numbers $\alpha^*$, $\beta^*$ such $\alpha^* = p^{*T}Aq^* + p^{*T}\pi$ and $\beta^* = - p^{*T}Aq^* + \rho^T q^*$.

Thus, $Aq^* + \pi \leq \alpha^* e^{(m)}$, $e^{(n)T}q^* = 1$, $q^* \in \mathbb{R}^n_+$, $\alpha^* \in \mathbb{R}$, and $-p^{*T}A + \rho^T \leq \beta^* e^{(n)T}$, $-p^{*T}e^{(m)} = -1$, $p^* \in \mathbb{R}^m_+$, $\beta^* \in \mathbb{R}$.

Further, $\rho^T q^* - \beta^* = p^{*T}Aq^* = \alpha^* - p^{*T}\pi$ and so $\rho^T q^* - \alpha^* = -\pi^T p^* + \beta^*$.

Suppose, $q^0$, $\alpha^0$ satisfy [$Aq^0 - \alpha^0 e^{(m)} \leq -\pi$, $e^{(n)T}q^0 = 1$, $q^0 \in \mathbb{R}^n_+$, $\alpha^0 \in \mathbb{R}$].

Since, [$-p^{*T}A + \rho^T \leq \beta^* e^{(n)T}$, $-p^{*T}e^{(m)} = -1$, $p^* \in \mathbb{R}^m_+$, $\beta^* \in \mathbb{R}$] from the first paragraph of this proof it follows that $\rho^T q^0 - \alpha^0 \leq -\pi^T p^* + \beta^* = \rho^T q^* - \alpha^*$.

Thus, $q^*$, $\alpha^*$ solve [Maximize $\rho^T q - \alpha$, subject to $Aq - \alpha e^{(m)} \leq -\pi$, $e^{(n)T}q = 1$, $q \in \mathbb{R}^n_+$, $\alpha \in \mathbb{R}$].

Now suppose, $p^0$, $\beta^0$ satisfy [$p^{0T}A + \beta^0 e^{(n)T} \geq \rho^T$, $-p^{0T}e^{(m)} = -1$, $p^0 \in \mathbb{R}^m_+$, $\beta^0 \in \mathbb{R}$].

Since, [$Aq^* + \pi \leq \alpha^* e^{(m)}$, $e^{(n)T}q^* = 1$, $q^* \in \mathbb{R}^n_+$, $\alpha^* \in \mathbb{R}$] from the first paragraph of this proof it follows that $-\pi^T p^0 + \beta^0 \geq \rho^T q^* - \alpha^* = -\pi^T p^* + \beta^*$.

Thus, $p^*$, $\beta^*$ solve [Minimize $-\pi^T p + \beta$, subject to $p^T A + \beta e^{(n)T} \geq \rho^T$, $-p^T e^{(m)} = -1$, $p \in \mathbb{R}^m_+$, $\beta \in \mathbb{R}$]. Further, $p^*$, $q^*$, $\alpha^*$, $\beta^*$ satisfies $\alpha^* = p^{*T}Aq^* + p^{*T}\pi$ and $\beta^* = -p^{*T}Aq^* + \rho^T q^*$.

**Part 2:** Now suppose $q^*$, $\alpha^*$ solve [Maximize $\rho^T q - \alpha$, subject to $Aq - \alpha e^{(m)} \leq -\pi$, $e^{(n)T}q = 1$, $q \in \mathbb{R}^n_+$, $\alpha \in \mathbb{R}$] and $p^*$, $\beta^*$ solve its dual [Minimize $-\pi^T p + \beta$, subject to $p^T A + \beta e^{(n)T} \geq \rho^T$, $-p^T e^{(m)} = -1$, $p \in \mathbb{R}^m_+$, $\beta \in \mathbb{R}$].

Then, we know that:

(1) $Aq^* - \alpha^* e^{(m)} \leq -\pi$, $e^{(n)T}q^* = 1$, $q^* \in \mathbb{R}^n_+$, $\alpha^* \in \mathbb{R}$ and $p^{*T}Aq^* - \alpha^* = -p^{*T}\pi$.

(2) $p^{*T}A + \beta^* e^{(n)T} \geq \rho^T$, $-p^{*T}e^{(m)} = -1$, $p^* \in \mathbb{R}^m_+$, $\beta^* \in \mathbb{R}$ and $p^{*T}Aq^* + \beta^* = \rho^T q^*$.

Thus,

(3) $\alpha^* = p^{*T}Aq^* + p^{*T}\pi$ and $\beta^* = -p^{*T}Aq^* + \rho^T q^*$.

Substituting for $\alpha^*$ and $\beta^*$ in (1) and (2) we get $Aq^* + \pi \leq (p^{*T}Aq^* + p^{*T}\pi)e^{(m)}$, $-p^{*T}A + \rho^T \leq (-p^{*T}Aq^* + \rho^T q^*)e^{(n)}$, $(p^*, q^*) \in \Delta^{m-1} \times \Delta^{n-1}$.

Thus, $(p^*, q^*) \in \Delta^{m-1} \times \Delta^{n-1}$ and for all $(p, q) \in \Delta^{m-1} \times \Delta^{n-1}$: $f^R(p^*, q^*) \geq f^R(p, q^*)$ and $f^C(p^*, q^*) \geq f^C(p^*, q)$.

Thus, $(p^*, q^*)$ is an equilibrium for the TPASS game $(A, \pi, \rho)$.

Further, from (4) we get $\alpha^* = p^{*T}Aq^* + p^{*T}\pi$ and $\beta^* = -p^{*T}Aq^* + \rho^T q^*$. Q.E.D.

Propositions 1 and 2 in this section is summarized in the following theorem.

**Theorem 1:** Let $(A, \pi, \rho)$ be a TPASS game.

(i) $(A, \pi, \rho)$ has an equilibrium.

(ii) $(p^*, q^*)$ is an equilibrium for the TPASS game $(A, \pi, \rho)$ <u>if and only if</u> there exist real numbers $\alpha^*$, $\beta^*$ such that $q^*$, $\alpha^*$ solve [Maximize $\rho^T q - \alpha$, subject to $Aq - \alpha e^{(m)} \leq -\pi$, $e^{(n)T}q = 1$, $q \in \mathbb{R}^n_+$, $\alpha \in \mathbb{R}$] and $p^*$, $\beta^*$ solve its dual [Minimize $-\pi^T p + \beta$, subject to $p^T A + \beta e^{(n)T} \geq \rho^T$, $-p^T e^{(m)} = -1$, $p \in \mathbb{R}^m_+$, $\beta \in \mathbb{R}$]. Further, $p^*$, $q^*$, $\alpha^*$, $\beta^*$ satisfies $\alpha^* = p^{*T}Aq^* + p^{*T}\pi$ and $\beta^* = -p^{*T}Aq^* + \rho^T q^*$.

**3.** The following result follows immediately from the "Equivalence Theorem" in section II of Mangasarian and Stone (1964). The Equivalence Theorem of Mangasarian and Stone, uses the result on existence of an equilibrium for a more general class of games, known as bi-matrix games, a proof of the latter being available as the proof of theorem 2 in Chandrasekaran (undated). Proposition 1 above, is a proof of existence for TPASS games that would hopefully be more widely accessible. It should be noted, that an existence result does not prove that the set of all equilibria for a TPASS game is equivalent to the set of all solutions of a linear programming problem. It is to the statement and proof of this latter result in the context of TPASS games that we turn to now.

**Proposition 3:** $(p^*, q^*)$ is an equilibrium for the TPASS game $(A, \pi, \rho)$ if and only if there exist real numbers $\alpha^*$, $\beta^*$ such that $p^*$, $q^*$, $\alpha^*$, $\beta^*$ solve the following "linear programming problem":

Maximize $\pi^T p + \rho^T q - \alpha - \beta$, subject to $Aq + \pi - \alpha e^{(m)} \leq 0$, $-A^T p + \rho - \beta e^{(n)} \leq 0$, $p^T e^{(m)} = 1$, $q^T e^{(n)} = 1$, $p \geq 0$, $q \geq 0$, $\alpha, \beta \in \mathbb{R}$.

**Proof:** First note that if $Aq + \pi - \alpha e^{(m)} \leq 0$, $-A^T p + \rho - \beta e^{(n)} \leq 0$, $p^T e^{(m)} = 1$, $q^T e^{(n)} = 1$, $p \geq 0$, $q \geq 0$, $\alpha, \beta \in \mathbb{R}$, then it must be the case that $\pi^T p + \rho^T q - \alpha - \beta \leq 0$.

This we get by pre-multiplying the first inequality by $p^T$, the second inequality by $q$ and adding the two resulting inequalities.

Suppose, $(p^*, q^*)$ be an equilibrium for the TPASS game $(A, \pi, \rho)$. Let $\alpha^* = p^{*T} A q^* + p^{*T} \pi$ and $\beta^* = -q^{*T} A^T p^* + q^{*T} \rho$.

Adding the two equalities we get $\pi^T p^* + \rho^T q^* - \alpha^* - \beta^* = 0$.

Since $(p^*, q^*)$ is an equilibrium it must be the case that $Aq^* + \pi - (p^{*T} A q^* + p^{*T} \pi) e^{(m)} \leq 0$, $-A^T p^* + \rho - (-p^* A q^* + q^{*T} \rho) e^{(n)} \leq 0$, $p^{*T} e^{(m)} = 1$, $q^{*T} e^{(n)} = 1$, $p^* \geq 0$, $q^* \geq 0$.

Thus, $p^*$, $q^*$, $\alpha^*$, $\beta^*$ satisfies all the constraints of the linear programming problem and $\pi^T p^* + \rho^T q^* - \alpha^* - \beta^* = 0$.

Thus, $p^*$, $q^*$, $\alpha^*$, $\beta^*$ solves the linear programming problem.

Now suppose $p^*$, $q^*$, $\alpha^*$, $\beta^*$ solves the linear programming problem.

From proposition 1, we know that the TPASS game has an equilibrium $(p^0, q^0)$ and we have just shown that since $(p^0, q^0)$ is an equilibrium, there exists real numbers $\alpha^0$, $\beta^0$ such that $p^0$, $q^0$, $\alpha^0$, $\beta^0$ solves the linear programming problem and $\pi^T p^0 + \rho^T q^0 - \alpha^0 - \beta^0 = 0$.

Thus, it must be the case that $\pi^T p^* + \rho^T q^* - \alpha^* - \beta^* = 0$.

Since, $Aq^* + \pi - \alpha^* e^{(m)} \leq 0$, $-A^T p^* + \rho - \beta^* e^{(n)} \leq 0$, $p^{*T} e^{(m)} = 1$, $q^{*T} e^{(n)} = 1$, $p^* \geq 0$, $q^* \geq 0$,

we get $p^{*T} A q^* + p^{*T} \pi - \alpha^* \leq 0$ and $-p^* A^T q^* + \rho^T q^* - \beta^* \leq 0$.

Adding the two inequalities, we get $\pi^T p^* + \rho^T q^* - \alpha^* - \beta^* \leq 0$.

Hence, $\pi^T p^* + \rho^T q^* - \alpha^* - \beta^* = 0$ implies $p^{*T} A q^* + p^{*T} \pi - \alpha^* = 0$ and $-p^* A^T q^* + \rho^T q^* - \beta^* = 0$.

Thus, $\alpha^* = p^{*T} A q^* + p^{*T} \pi$ and $\beta^* = -p^* A^T q^* + \rho^T q^*$.

Thus, $Aq^* + \pi - (p^{*T} A q^* + p^{*T} \pi) e^{(m)} \leq 0$, $-A^T p^* + \rho - (-p^* A^T q^* + \rho^T q^*) e^{(n)} \leq 0$.

Thus, $(p^*, q^*)$ is an equilibrium for the TPASS game. Q.E.D.

**Note:** Chakrabarti, Gilles and Mallozzi (2024) introduce two classes of non-cooperative games in normal form with the number of players being finite but at least two. The first is the class of “separable games” and the second is the class of “additively separable games”. Our TPASS games are “a special case” in the two-person (two-player) context of “additively separable games”, the latter being definition 4.1 of Chakrabarti, Gilles and Mallozzi (2024). Since, the two-person version of games considered in Chakrabarti, Gilles and Mallozzi (2024)- including their ““additively separable games”- may not be bi-matrix games, the possibility of extending either Proposition 2 or Proposition 3 in our paper to the class of two-person “additively separable” games is an open question.

**4.** We now provide an example of a 2×2 TPASS game that “resembles” the well-known “Prisoners’ Dilemma”.

Consider the TPASS game (A, π, ρ), with A = $\begin{bmatrix} 0 & 1 \\ -1 & 0 \end{bmatrix}$, π, ρ$\in \mathbb{R}^2_{++}$ such that π = ρ and 1 > $\pi_2 > \pi_1 > 0$. For example, let $\pi_2 = \frac{3}{4}$, $\pi_1 = \frac{1}{2}$.The pay-off matrix for the “row player” is $\begin{bmatrix} \frac{1}{2} & \frac{3}{2} \\ -\frac{1}{4} & \frac{3}{4} \end{bmatrix}$ and the pay-off matrix for the “column player” is $\begin{bmatrix} \frac{1}{2} & -\frac{1}{4} \\ \frac{3}{2} & \frac{3}{4} \end{bmatrix}$.

In this case if we interpret the first strategies of both players as “non-cooperation” and the second strategy of both players as “cooperation”, then the pay-off matrices for the two players in the TPASS game (A, π, ρ) resemble a “Prisoner’s Dilemma” situation.

**Acknowledgment:** I wish to thank Subhadip Chakrabarti for drawing my attention to his and his co-authors’ work on additively separable games. I would also like to thank Bernhard von Stengel for a “cautionary” example related to an earlier version of the paper that has now been taken into consideration.